\documentclass[10pt, a4paper, reqno]{amsart}
\usepackage{preambolo}

\begin{document}


\setcounter{secnumdepth}{3}

\setcounter{tocdepth}{2}

\title{\textbf{Invariant and non-invariant almost complex structures on compact quotients of Lie groups}}

\author[Lorenzo Sillari]{Lorenzo Sillari}

\address{Lorenzo Sillari: Dipartimento di Scienze Matematiche, Fisiche e Informatiche, Unità di Matematica e Informatica, Università degli Studi di Parma, Parco Area delle Scienze 53/A, 43124, Parma, Italy} \email{lorenzo.sillari@unipr.it}

\author[Adriano Tomassini]{Adriano Tomassini}

\address{Adriano Tomassini: Dipartimento di Scienze Matematiche, Fisiche e Informatiche, Unità di Matematica e Informatica, Università degli Studi di Parma, Parco Area delle Scienze 53/A, 43124, Parma, Italy}
\email{adriano.tomassini@unipr.it}

\maketitle

\begin{abstract} 
\noindent \textsc{Abstract}. In this paper we briefly survey the classical problem of understanding which Lie algebras admit a complex structure, put in the broader perspective of almost complex structures with special properties. We focus on the different behavior of invariant and non-invariant structures, with a special attention to their canonical bundle and Kodaira dimension. We provide new examples of computations of Kodaira dimension of invariant and non-invariant structures.
\end{abstract}

\blfootnote{  \hspace{-0.55cm} 
{\scriptsize 2020 \textit{Mathematics Subject Classification}. Primary: 32Q60, 53C15; Secondary: 22E25, 22E60. \\ 
\textit{Keywords: canonical bundle, invariant almost complex structure, invariant cohomology, Iwasawa manifold, Kodaira dimension, Nakamura manifold, nilpotent Lie algebra, rank of the Nijenhuis tensor, solvable Lie algebra.}\\}

\noindent The authors are partially supported by GNSAGA of INdAM. The second author is partially supported by the Project PRIN 2022 “Real and Complex Manifolds: Geometry and Holomorphic Dynamics” (code 2022AP8HZ9)}

\section{Introduction}\label{sec:intro}

In \cite{Has05}, Hasegawa classifies complex surfaces whose underlying manifold is diffeomeorphic to a $4$-dimensional solvmanifold. He also shows that every complex structure on such manifolds is biholomorphic to an invariant one, and he points out that, in higher dimension, there is no known example of a solvmanifold admitting a complex structure that is not biholomorphic to an invariant one. Hence, he poses the natural
\vspace{.2cm}

\textbf{Question:} is every complex structure on a solvmanifold biholomorphic to an invariant one?
\vspace{.2cm}

Some years later, Hasegawa provides a negative answer: he shows that every invariant complex structure on the Nakamura manifold $\N$, a $6$-dimensional solvmanifold, is biholomorphic to $\C^3$ \cite{Has10}, while Nakamura shows that $\N$ admits (non-invariant) complex structures whose universal cover is not Stein \cite{Nak75}. In particular, they cannot be biholomorphic to the invariant ones. The same conclusion can be reached thinking in terms of the Kodaira dimension of complex structures, which is invariant under biholomorphisms. Since every invariant complex structure on $\N$ has $\C^3$ as a universal cover, the holomorphic $(3,0)$-form $dz^1 \wedge dz^2 \wedge dz^3$ induces a holomorphic trivialization of the canonical bundle of the quotient, implying that every invariant structure has Kodaira dimension $0$. On the other side, Kodaira (see \cite{Nak75}), built a non-invariant deformation of the standard complex structure on $\N$ with Kodaira dimension $- \infty$.

The example we described shows that invariant and non-invariant structures on compact quotients of Lie groups might exhibit very different behaviors (see also Hasegawa's survey \cite{Has09}). This paper aims to be a short survey on invariant complex structures, their cohomology, and their Kodaira dimension, within the broader perspective of invariant \textit{almost} complex structures, and to show with explicit examples that the difference between the invariant and the non-invariant case is a phenomenon that is not limited to integrable structures. For instance, the rank of the Nijenhuis tensor of invariant structures has an upper bound given in terms of the Lie algebra cohomology, while this is not the case for the rank of non-invariant structures, see Theorem 5.15 in \cite{ST23b} and the examples therein.

The paper is structured as follows. In Section \ref{sec:prel}, we recall some basics on invariant almost complex structures, on compact quotients of Lie groups, and on the almost complex Kodaira dimension. Section \ref{sec:acs:algebra} contains a fast review of the literature on invariant complex structures and some recent progress in the almost complex case. Finally, the examples of Section \ref{sec:noninv} are new material and provide explicit computations of invariant and non-invariant structures on the Iwasawa manifold and of their Kodaira dimension.

\section{Preliminaries}\label{sec:prel}

\subsection*{Notation.}

Given a basis of complex forms $\phi^j$, $j=1, \ldots, m$, we abbreviate $\phi^j \wedge \phi^k$ to $\phi^{jk}$ and $\bar \phi^j$ to $\phi^{\bar j}$.

\subsection*{Lie algebras and almost complex structures.}

Let $\g = \R \langle e_1, \ldots, e_{2m} \rangle$ be a $2m$-dimensional real Lie algebra and let $\g^* = \R \langle e^1, \ldots, e^{2m} \rangle$ be its dual. Consider the graded algebra
\[
A^k_\g \coloneqq \bigwedge^k \g^*
\]
endowed with the differential defined as the dual of the Lie bracket: for $X$, $Y \in \g$ and $a \in \g^*$, we set
\[
da (X,Y) \coloneq - a ( [X,Y] ),
\]
and extend it as a graded derivation to $A^\bullet_\g$. The differential graded algebra $(A^\bullet_\g, d)$ is the \emph{Chevalley--Eilenberg complex of $\g$} \cite{CE48}. 

Let $J \in \End (\g)$, with $J^2 = - \Id$, be an almost complex structure on $\g$. Elements of the complexified Chevalley--Eilenberg complex $A^k_{\g^\C} \coloneqq \bigwedge^k (\g^*)^\C$ inherit a bigrading induced by $J$ via the decomposition $\g^\C = \g^{1,0} \oplus \g^{0,1}$ into $(\pm i)$-eigenspaces. With respect to the bigrading, the differential splits as $d = \mu + \partial + \bar \partial + \bar \mu$, and the almost complex structure $J$ is \textit{integrable}, that is, its Nijenhuis tensor
\[
N_J (X,Y) \coloneqq [JX,JY] - J[JX,Y] - J[X,JY]-[X,Y], \quad X, Y \in \g,
\]
vanishes, if and only if $\bar \mu =0$. In this case, we have that $d = \partial + \bar \partial$ and we call $J$ a \emph{complex structure} on $\g$.

Complex structures satisfying nicer algebraic properties are often considered in the literature. In particular, $J$ is called \textit{abelian} if $[ \g^{1,0}, \g^{1,0} ] =0$, while it is called \textit{nilpotent} if $\mathfrak{a}_k (J) = \g$ for some $k$, where $\mathfrak{a}_k (J)$ is the $J$-compatible ascending series given by
\[
\mathfrak{a}_0 = \{ 0 \} \quad \text{and} \quad \mathfrak{a}_{k+1} \coloneqq \{ X \in \g: [X, \g] \subseteq \mathfrak{a}_k \text{ and } [JX, \g] \subseteq \mathfrak{a}_k \} \text{ for $k \ge 0$}.
\]

\subsection*{Structures on compact quotients of Lie groups.} 

Let $G$ be a connected, simply connected Lie group. Suppose that $G$ admits a co-compact lattice $\Gamma$ and denote by $M \coloneqq \Gamma \backslash G$ the compact quotient. If $G$ is nilpotent, resp.\ solvable, we say that $M$ is a \emph{nilmanifold}, resp. a \textit{solvmanifold}.

Let $J$ be an almost complex structure on $G$. The almost complex structure $J$ is \emph{left-invariant} if $L_h^* J = J L_h^*$ for all $h \in G$, where $L_h^*$ is the pullback by left multiplication. If $L_\gamma^* J = J L_\gamma^*$ for all $\gamma \in \Gamma$, we say that $J$ is \emph{$\Gamma$-invariant}. 

A smooth form $\alpha$ on $G$ is said to be \emph{left-invariant}, resp.\ \emph{$\Gamma$-invariant}, if $L_h^* \alpha = \alpha$ for all $h \in G$, resp.\ for all $h \in \Gamma$. Left-invariant (real) forms can be canonically identified with elements of the Chevalley--Eilenberg complex $A^\bullet_\g$, while left-invariant almost complex structures can be identified with almost complex structures on $\g$. The same holds for complex forms and the complex $A^k_{\g^\C}$. In addition, functions and structures on the quotient $M = \Gamma \backslash G$ have a precise counterpart in terms of functions and structures on $G$:

\begin{itemize}

    \item the space $C^\infty (M)$ can be identified with smooth functions on $G$ that are also $\Gamma$-periodic;

    \item almost complex structures on $M$ that are induced by left-invariant structures on $G$ are called \emph{invariant structures}. They correspond to almost complex structures on $\g$ or, equivalently, to linear maps $J$ sending each $e_j$ to $\R \langle e_1, \ldots, e_{2m} \rangle$ such that $J^2 = -\Id$;
    
    \item general almost complex structures on $M$ can be thought as $\Gamma$-invariant structures on $G$ or, equivalently, as $C^\infty(M)$-linear maps $J$ sending each $e_j$ to $C^\infty (M) \langle e_1, \ldots, e_{2m} \rangle$ such that $J^2 = -\Id$.
\end{itemize}

Due to the correspondence between almost complex structures on $\g$ and invariant almost complex structures on $M$, when dealing with invariant structures we will always work with structures on $\g$. 

\subsection*{The almost complex Kodaira dimension.}

Let $(M,J)$ be a compact almost complex manifold. The \emph{canonical bundle} of $J$ is the line bundle of $(m,0)$-forms $K_J \coloneqq \Lambda^{m,0}$. The space of its smooth sections is the space of $(m,0)$-forms. A section $\Omega$ of $K_J$ is said to be \emph{pseudoholomorphic} if $\bar \partial \Omega =0$. More in general, we call pseudoholomorphic any $\bar \partial$-closed form. If $K_J$ admits a never-vanishing section $\Omega$, then it is trivial as a complex line bundle. If, in addition, we can take $\Omega$ to be pseudoholomorphic, then $K_J$ is \textit{pseudoholomorphically trivial}, cf.\ \cite{BT96}. The choice of a (not necessarily pseudoholomorphic) trivialization of $K_J$ provides an $\mathrm{SU}(m)$-structure on $M$.

Denote by $P_1$ the space of pseudoholomorphic smooth sections of $K_J$, i.e., the space of pseudoholomorphic $(m,0)$-forms, and by $P_l$ the space of pseudoholomorphic sections of $K_J^{\otimes l}$, for $l \ge 1$. The \emph{Kodaira dimension of an almost complex manifold}, see \cite{CZ23} and \cite{CZ24}, is the number
    \[
    \kappa_J := 
    \begin{cases}
        - \infty &\text{if $P_l = \{ 0 \}$ for all $l \ge 1$,}\\
        \limsup\limits_{l \to \infty} \frac{\log \dim_\C P_l}{\log l} &\text{otherwise}.
    \end{cases}
    \]
In particular, observe that if the canonical bundle is pseudoholomorphically trivial, then the Kodaira dimension of $J$ vanishes \cite{CNT24}.

\section{Almost complex structures on Lie algebras}\label{sec:acs:algebra}
 
In this section we review some problems related to almost complex structures on Lie algebras. We focus on three aspects: existence, classification, and cohomologies. Due to the importance that the Kodaira dimension plays in complex and almost complex geometry, we also briefly review invariant complex structures with holomorphically trivial canonical bundle.

\subsection{Existence and classification of complex structures.}\label{sec:complex}

Let $\g$ be a $2m$-dimensional Lie algebra and let $J$ be a complex structure on $\g$. The two main problems when dealing with complex structures on Lie algebras are the existence, i.e., understanding which Lie algebras admit complex structures, and the classification, i.e., determine the (holomorphic) isomorphism classes of complex structures on a given Lie algebra.

These questions are hard to handle in full generality. To simplify the approach, one usually makes further assumptions on $\g$ or $J$. We give an overview of the literature on the existence and classification problem for complex structures on Lie algebras. 

We first point out the survey by Smolentsev, which collects results in dimensions $4$ and $6$, together with detailed references to classical results on the existence problem \cite{Smo15}. 

The most studied case is that of nilpotent Lie algebras. In dimension $4$, it is an exercise to determine which nilpotent Lie algebras admit a complex structure. In dimension $6$, there are 34 isomorphism classes of nilpotent Lie algebras. Salamon shows that $18$ of them admit a complex structure \cite{Sal01}. Together with Ketsetzis, he fully classifies invariant complex structures on the Iwasawa manifold up to isomorphism \cite{KS04}, while Magnin does so on indecomposable $6$-dimensional nilpotent Lie algebras \cite{Mag07}, and Ceballos, Otal, Ugarte and Villacampa on arbitrary $6$-dimensional nilpotent Lie algebras \cite{COUV16}. General results in higher dimension are scarce. For dimension $8$, we point out some recent papers by Latorre, Ugarte and Villacampa on arbitrary nilpotent Lie algebras \cite{LUV19}, and on nilpotent Lie algebras with one-dimensional center \cite{LUV23}, and by Millionshchikov on nilpotent Lie algebras with $b_1 (\g) = 2$ \cite{Mil24}. The last two papers include a full classification.

In the solvable case, the problems of existence and classification are solved, in dimension $4$, by Snow \cite{Sno90} and Ovando \cite{Ova00}. There are no general results in higher dimension.

Regarding special complex structures, Cordero, Fern\'andez, Gray and Ugarte deal with nilpotent complex structures on nilmanifolds, see \cite{CFGU00}, \cite{CFGU01} and \cite{CFGU97}. Angella, Otal, Ugarte and Villacampa study complex structures of splitting type on $6$-dimensional solvable Lie algebras \cite{AOUV17}. Special attention is dedicated to the study of abelian complex structures by Cordero, Fern\'andez and Gray on $6$-dimensional nilpotent Lie algebras \cite{CFU02}, by Barberis and Dotti on solvable Lie algebras \cite{BD04}, and by Andrada, Barberis and Dotti on $6$-dimensional Lie algebras with a full classification \cite{ABD11}, see also \cite{ABD13}.

In addition to the nilpotent and solvable cases, one can find results on other special classes of Lie algebras. Goze and Remm focus on filiform Lie algebras \cite{GR02}, Garc\'ia Vergnolle and Remm on quasi-filiform Lie algebras \cite{GR09}, Czarnecki and Sroka on $6$-dimesnional Lie algebras that are a product \cite{CS18}, and Sroka on $6$-dimensional decomposable Lie algebras \cite{Sro20}.

\subsection{Existence and classification of almost complex structures with special properties.}\label{sec:acs}

Every $2m$-dimensional Lie algebra $\g$ admits many almost complex structures, which are parametrized by $GL (2m, \R) / GL(m,\C)$. Hence, to obtain meaningful existence and classification results, one has to impose further conditions on $J$.

Let $J$ be an almost complex structure on $\g$ and let $N_J$ be its Nijenhuis tensor. Since $N_J$ satisfies the symmetry
\[
N_J(J X, Y) = N_J (X, JY) = - J N_J(X,Y), \quad X,Y \in \g,
\]
its image is an even-dimensional vector space $\mathcal{V}$, whose complex dimension defines an invariant of $J$
\[
\rk N_J \coloneqq \dim_\C \mathcal{V},
\]
called the \emph{rank of the Nijenhuis tensor}. The most natural problem is to determine which Lie algebras admit structures of a given rank. In particular, observe that almost complex structures with $\rk N_J =0$ correspond to complex structures, that are discussed in detail in Section \ref{sec:complex}. The literature on the almost complex case is scarcer than the one on integrable structures. Mushkarov classifies nilpotent Lie algebras of dimension $4$ and $6$ for which every almost complex structure has at least one non-constant pseudoholomorphic function \cite{Mus14}. We also point out that in \cite{ST23b}, the authors expand Salamon's classification by determining all the possible ranks of almost complex structures on $6$-dimensional nilpotent Lie algebras. The classification shows a very different behavior between invariant and non-invariant structures in terms of the possible values of the rank. For instance, on a torus every invariant structure is integrable, while there are plenty of non-invariant structures that are not integrable. On the other side, on the Iwasawa manifold, every invariant structure $J$ satisfies, by Theorem 5.15 in \cite{ST23b}, the topological bound $\rk N_J \le 6 - b_1 = 2$, while there are non-invariant structures with $\rk N_J=3$.
We also point out an example, due to Fern\'andez, de Le\'on and Saralegui, of a $6$-dimensional solvmanifold that does not admit invariant complex structures \cite{FLS96}. It is not known if it admits non-invariant ones. 

Under additional assumptions, it is possible to classify almost complex structures up to isomorphism. We set
\[
\mathcal{V}^{(0)} \coloneqq \mathcal{V} \quad \text{and} \quad \mathcal{V}^{(k+1)} \coloneqq \mathcal{V}^{(k)} + [ \mathcal{V}^{(k)}, \mathcal{V}^{(k)}], \quad \text{ for $k \ge 0$}.
\]
We say that $\mathcal{V}$ is \emph{fundamental} if $\mathcal{V}^{(k)} = \g$ for some $k$. In such a case, the almost complex structure $J$ admits a canonical form in terms of a basis of $\g$ \cite{Kru98}. Kruglikov explicitly provides the canonical form in dimensions $4$ and $6$ \cite{Kru14}, while Bozzetti and Medori classify the structures, up to isomorphism, on $4$-dimensional non-solvable Lie algebras \cite{BM17}.

\subsection{Invariant cohomologies.}\label{sec:cohomology}

Let $M= \Gamma \backslash G$ be a compact quotient of a Lie group by a lattice. Denote by $H^k_d$ the de Rham cohomology of $M$, and, if $M$ admits an invariant complex structure, denote by $H^{p,q}_{\bar \partial}$ its Dolbeault cohomology. 

The \emph{invariant de Rham cohomology}, resp.\ \textit{invariant Dolbeault cohomology}, of $M$ is the de Rham cohomology, resp.\ Dolbeault cohomology, computed on invariant forms. We denote it by $\prescript{L}{}{H}^k_d$, resp.\ by $\prescript{L}{}{H}^{p,q}_{\bar \partial}$.

In general, there are injections of invariant cohomologies into the genuine cohomologies
\[
\prescript{L}{}{H}^k_d \longhookrightarrow {H}^k_d \quad \text{and} \quad \prescript{L}{}{H}^{p,q}_{\bar \partial} \longhookrightarrow {H}^{p,q}_{\bar \partial},
\]
and understanding if the maps are isomorphisms, i.e., if the de Rham and Dolbeault cohomologies can be computed using only invariant forms, is the object of open problems. Regarding the de Rham cohomology, Nomizu showed that the isomorphism holds for nilmanifolds \cite{Nom54}, while Hattori showed that it holds for solvmanifolds whose corresponding Lie algebra is completely solvable \cite{Hat60}. De Bartolomeis and the the second author built an example of a non-completely solvable solvmanifold where the isomorphism does not hold \cite{BT06a}, see also the work of Console and Fino \cite{CF11}.

Regarding the Dolbeault cohomology, it is known that, for arbitrary solvmanifolds, the inclusion of invariant Dolbeault cohomology into the Dolbeault cohomology is not an isomorphism. For instance, the Nakamura manifold admits a holomorphic parallelism given by invariant holomorphic $(1,0)$-forms $\phi^j$, $j=1,2,3$. Its Dolbeault cohomology group $H^{0,1}_{\bar \partial}$ has complex dimension $3$ and it is generated by $\phi^1$, which is invariant, and by two non-invariant forms obtained by multiplying $\phi^2$ and $\phi^3$ by smooth functions \cite{Nak75}. The isomorphism fails also in the completely solvable case, as shown by Kasuya on the completely solvable Nakamura manifold \cite{Kas13}.

In the nilpotent case, it is conjectured that the inclusion is always an isomorphism. It seems that the conjecture first appeared as a question in the works of Console and Fino \cite{CF01} and of Cordero, Fern\'andez, Gray and Ugarte \cite{CFGU00}. The conjecture is discussed in great detail by Rollenske \cite{Rol11}, and it is known to hold in several cases, see Theorem 3.2 in \cite{Ang14}, \cite{FRR19}, \cite{Rol11}, \cite{Rol09} and the references therein. In particular, it holds for nilmanifolds up to dimension $6$.

In the almost complex case, there are several cohomologies defined by Cirici and Wilson \cite{CW21}, by Coelho, Placini and Stelzig \cite{CPS22}, by the authors, see \cite{ST23a} and \cite{ST24}, and by Cahen, Gutt and Gutt \cite{CGG24}, that one can consider. Unfortunately, a general isomorphism between invariant and non-invariant cohomologies cannot hold in these cases, since the former is always finite-dimensional, while the latter is infinite-dimensional even on almost complex tori \cite{CPS22}.

\subsection{The canonical bundle of invariant almost complex structures.}\label{sec:canonical}

Let $\g$ be a $2m$-dimensional Lie algebra, let $J$ be an almost complex structure on $\g$ and let $K_J = \bigwedge^m \g^{1,0}$ be its canonical bundle. In contrast to what could happen for non-invariant structures on compact quotients of Lie group, if $J$ is invariant then $K_J$ is necessarily trivial since, if $\phi^j$, $j=1, \ldots, m$, is a basis of $\g^{1,0}$, then $\Omega \coloneqq \phi^{1 \ldots m}$ is a never-vanishing section of $K_J$. In general, we have that $\bar \partial \Omega \neq 0$ and $K_J$ is not pseudoholomorphically trivial.

Barberis, Dotti and Verbitsky show that every \textit{complex} structure on a nilpotent Lie algebra has holomorphically trivial canonical bundle \cite{BDV09}, see also \cite{CG04}, and this is the case also for certain solvable Lie algebras, see the examples in dimension $6$ by Fino, Otal and Ugarte \cite{FOU15} for the complex case, and the examples in dimension $4$ by Cattaneo \cite{Cat22}, and by Cattaneo, Nannicini and the second author \cite{CNT24}, \cite{CNT20} and \cite{CNT21} for the almost complex case. A full classification of $6$-dimensional unimodular Lie algebras admitting a complex structure with holomorphically trivial canonical bundle is given by Otal and Ugarte \cite{OU23}.

It is important to remark that, when working on compact quotients of Lie groups, there are invariant (almost) complex structures $J$ that do not admit an invariant (pseudo)holomorphic trivialization, but admit a non-invariant one, see Example \ref{ex:trivial} for some explicit computation, or \cite{AT24} and \cite{Tol24} for general results.

On the almost complex side, we recall that pseudoholomorphic triviality of the canonical bundle has been considered before by de Bartolomeis and Tian in the setting of \emph{bundle almost complex structures} \cite{BT96}, and by de Bartolomeis and the second author in the context of \emph{QIS manifolds}, see \cite{BT13} and \cite{BT06b}.

\section{Non-invariant structures}\label{sec:noninv}

In this section we build explicit families of non-invariant complex and almost complex structures on the Iwasawa manifold.

\subsection{Non-invariant complex structures on the Iwasawa manifold with Kodaira dimension \texorpdfstring{$0$}{}.}

Let $\mathbb{H}_3^\C$ be the complex Heisenberg group
\[
\mathbb{H}_3^\C \coloneqq \left \{ 
\begin{bmatrix}
1 & z_1 & z_3 \\
  & 1 & z_2 \\
  &   & 1 \\
\end{bmatrix}
: z_1,z_2,z_3 \in \C
\right\}
\]
whit the group operation induced by matrix multiplication. The Iwasawa manifold is the $6$-dimensional nilmanifold obtained as the quotient
\[
M \coloneqq ( \mathbb{H}_3^\C \cap {\mathrm {SL}} (3, \Z[i]) ) \backslash \mathbb{H}_3^\C.
\]
The complex structure inherited from $\C^3$ induces on $M$ a basis of $(1,0)$-forms $\phi^j$, $j=1, 2,3$, whose differentials are
\[
d \phi^1 =0, \quad d \phi^2=0 \quad \text{and} \quad d \phi^3 = - \phi^{12}.
\]
A basis of $(1,0)$-vector fields is given by
\begin{equation}\label{eq:xi}
\xi_1 = \partial_{z_1}, \quad \xi_2 = \partial_{z_2} + z_1 \partial_{z_3} \quad \text{and} \quad \xi_3 = \partial_{z_3}.
\end{equation}
Let $A \in \R$, $\abs{A} < 1$, and consider the family of almost complex structures $J_A$ defined by the coframe of $(1,0)$-forms
\[
\omega^1 \coloneqq \phi^1, \quad \omega^2 \coloneqq \phi^2 + A \sin (x_2) \phi^{\bar 2} \quad \text{and} \quad \omega^3 \coloneqq \phi^3,
\]
where $x_2 = \Re z_2$. Set $B \coloneqq (1- \abs{A}^2 \abs{\sin (x_2)}^2)^{-1}$. The frame of $(1,0)$-vector fields dual to the $\omega^j$ is given by 
\[
\psi_1 = \xi_1, \quad \psi_2 = B (\xi_2 - A \sin (x_2) \xi_{\bar 2}) \quad \text{and} \quad \psi_3 = \xi_3.
\]
We have the differentials
\[
d \omega^1 =0, \quad d \omega^2 = A B \, \partial_{z_2} ( \sin (x_2)) \omega^{2 \bar 2} \quad \text{and} \quad d\omega^3 = - B (\omega^{12} - A \sin(x_2) \omega^{1 \bar 2}),
\]
so that the structure $J_A$ is integrable. It is immediate to see that
\[
\bar \partial (\omega^{123}) = B \, \partial_{z_2} (\sin (x_2)) \omega^{123 \bar 2},
\]
and that the canonical bundle of $J_A$ is not holomorphically trivialized by $\omega^{123}$. Nevertheless, an easy check shows that $f = (1 - A \sin (x_2))^{-1}$ is a solution of the equation
\[
\bar \partial (f \omega^{123} ) = 0,
\]
and that the canonical bundle of $J_A$ is holomorphically trivial. More in general, if $p \in C^\infty (M)$ is a non-constant real function depending only on the $x_2$ variable and such that $\abs{p} \neq 1$ on $M$, we can consider the structure $J_{p}$ defined by 
\[
\omega^1 \coloneqq \phi^1, \quad \omega^2 \coloneqq \phi^2 + p \, \phi^{\bar 2} \quad \text{and} \quad \omega^3 \coloneqq \phi^3,
\]
Then the form $ (1-p)^{-1} \omega^{123}$ provides a holomorphic trivialization of the canonical bundle of $J_{p}$. In this case, the Kodaira dimension is $\kappa_{J_p} = 0$. We observe that this is coherent with Conjecture (ii) in \cite{Has10}, stating that every small deformation of an invariant complex structure on a nilmanifold is biholomorphic to an invariant one.

\subsection{Invariant and non-invariant almost complex structures on the Iwasawa manifold with Kodaira dimension \texorpdfstring{$0$}{} or \texorpdfstring{$- \infty$}{}.} We give three examples of almost complex structures on the Iwasawa manifold. The first one is a family of invariant structures $J_{A,B}$, $A, B \in \C$, whose canonical bundle is pseudoholomorphically trivialized by an invariant form, hence $\kappa_{J_{A,B}} = 0$. The second one is an invariant structure $J_0$ whose canonical bundle is pseudoholomorphically trivialized by a non-invariant form, but not by an invariant one, giving $\kappa_{J_0}=0$. The last one is a family of non-invariant structures $J_{q,r}$, $q,r \in C^\infty (M)$, whose canonical bundle is not pseudoholomorphically trivial, so that $\kappa_{J_{q,r}} = -\infty$.

Despite several attempts, we have not been able to produce an \textit{invariant} almost complex structure on the Iwasawa manifold with Kodaira dimension $- \infty$. Thus, we ask the following.
\vspace{.2cm}

\textbf{Question:} does every invariant almost complex structure on the Iwasawa manifold (or, more generally, on a nilmanifold) necessarily have pseudoholomorphically trivial canonical bundle?

\begin{example}
    Let $J_{A,B}$, with $A,B \in \C$, $\abs{A} \neq 1$ and $\abs{B} \neq 1$, be the family of invariant almost complex structures on the Iwasawa manifold defined by the co-frame
    \[
    \omega^1 \coloneqq \phi^1 + A \phi^{\bar 1}, \quad \omega^2  \coloneqq  \phi^2 + B \phi^{\bar 2} \quad \text{and} \quad \omega^3  \coloneqq  \phi^3.
    \]
    The differentials are
    \begin{align*}
    &d \omega^1 =0, \quad d \omega^2 =0 \quad \text{and}\\
    &d \omega^3 = - (1-\abs{A}^2)^{-1}(1-\abs{B}^2)^{-1} (\omega^{12}- B \omega^{1 \bar 2} - A \omega^{\bar 1 2} + AB \omega^{\bar 1 \bar 2} ).
    \end{align*}
    It is immediate to see that $J_{A,B}$ is not integrable and that $\bar \partial \omega^{123} =0$. \hfill $\blacksquare$
\end{example}

\begin{example}\label{ex:trivial}
    Let $e^j$, $j=1, \ldots, 6$, be the standard real basis for the dual of the Lie algebra associated to the Iwasawa manifold, with structure equations
    \[
    de^j =0, \quad \text{for $j =1, \ldots, 4$}, \quad de^5 = e^{13} - e ^{24} \quad \text{and} \quad de^6 = e^{14} + e ^{23}.
    \]
    Consider the almost complex structure $J_0$ defined by the co-frame
    \[
    \omega^1 \coloneqq e^1 + i e^2, \quad \omega^2  \coloneqq  e^3 + i e^5 \quad \text{and} \quad \omega^3  \coloneqq  e^4 + i e^6.
    \]
    The differentials are
    \begin{align*}
    &d \omega^1 = 0, \quad d\omega^2 = - \frac{1}{4} (\phi^{13} + \phi^{1 \bar 3}- \phi^{\bar 1 3}- \phi^{\bar 1 \bar 3}) + \frac{i}{4} (\phi^{12} + \phi^{1 \bar 2}+ \phi^{\bar 1 2} + \phi^{\bar 1 \bar 2})\\
    &\text{and} \quad d\omega^3 = \frac{1}{4} (\phi^{12} + \phi^{1 \bar 2}- \phi^{\bar 1 2}- \phi^{\bar 1 \bar 2}) + \frac{i}{4} (\phi^{13} + \phi^{1 \bar 3}+ \phi^{\bar 1 3} + \phi^{\bar 1 \bar 3}).
    \end{align*}
    Then $J_0$ is not integrable and 
    \[
    \bar \partial \omega^{123} = - \frac{i}{2} \phi^{123\bar 1},
    \]
    so that its canonical bundle is not trivialized by an invariant form. However, the non-invariant form $e^{ix_1} \omega^{123}$ provides such a trivialization, where $x_1 = \Re z_1 $. \hfill $\blacksquare$
\end{example}

\begin{example}
    Let $J_{q,r}$, with $q,r \in C^\infty(M)$ to be chosen later, be the family of non-invariant almost complex structures on the Iwasawa manifold defined by the co-frame
    \[
    \omega^1 \coloneqq \phi^1, \quad \omega^2  \coloneqq  \phi^2 \quad \text{and} \quad \omega^3  \coloneqq  \phi^3 + q \phi^{\bar 1} + r \phi^{\bar 2}.
    \]
    The differentials are
    \[
    d \omega^1 =0, \quad d \omega^2 =0 \quad \text{and} \quad d \omega^3 = - \omega^{12} + dq \wedge \omega^{\bar 1} + dr \wedge \omega^{\bar 2}.
    \]
    Let $\xi_j$, $j =1,2,3$, resp.\ $\psi_j$, $j = 1,2,3$, be co-frames of $(1,0)$-vector fields for the $\phi^j$, resp.\ the $\omega^j$. Then, we have that
    \[
    \psi_1 = \xi_1 - \bar q \, \xi_{\bar 3}, \quad \psi_2 = \xi_2 - \bar r \, \xi_{\bar 3} \quad \text{and} \quad \psi_3 = \xi_3, 
    \]
    where the explicit expression of the $\xi_j$ is given by \eqref{eq:xi}. To show that the canonical bundle of $J_{q,r}$ is not pseudoholomorphically trivial, we need to show that the equation $\bar \partial (f (\omega^{123})^{\otimes l}) =0$ has no solution $f \in C^\infty (M)$ for all $l \in \N$. In the case $l =1$, we impose
    \[
    \bar \partial (f \omega^{123})  = \bar \partial f \wedge \omega^{123} + f \omega^{12} \wedge \bar \partial \omega^3 = \bar \partial f \wedge \omega^{123} +f \omega^{12} \wedge ( \partial q \wedge \omega^{\bar 1} + \partial r \wedge \omega^{\bar{2}}) = 0,
    \]
    which gives the system of equations
    \[
    \begin{cases}
        \psi_{\bar 1} (f) = \psi_3 (q) f, \\
        \psi_{\bar 2} (f) = \psi_3 (r) f, \\
        \psi_{\bar 3} (f) = 0.
    \end{cases}
    \]
    Since $\psi_{\bar 3} = \xi_{\bar 3} = \partial_{\bar z_3}$, from the last equation we get that $\partial_{z_3} \partial_{\bar z_3} f =0$. The projection on the $z_1$ and $z_2$ coordinates gives a fibration of the Iwasawa manifold as a $T^2$-bundle over $T^4$. On the fiber of the fibration, the operator $\partial_{z_3} \partial_{\bar z_3}$ is elliptic and, by compactness of the fibers, $f$ must be constant on each fiber. Thus, $f$ does not depend on the $z_3$ variable, and the system reduces to
    \[
    \begin{cases}
        \xi_{\bar 1} (f) = \xi_3 (q) f, \\
        \xi_{\bar 2} (f) = \xi_3 (r) f.
    \end{cases}
    \]
    It is enough to choose $q$ or $r$ \textit{not independent} of the $z_3$ variable to conclude that the system does not admit solutions $f \in C^\infty (M)$. In a similar way, a function $f$ is a solution of the equation $\bar \partial ( f (\omega^{123})^{\otimes l}) =0$ if and only if it is independent of the $z_3$ variable and it solves the system
    \[
    \begin{cases}
        \xi_{\bar 1} (f) = l \xi_3 (q) f, \\
        \xi_{\bar 2} (f) = l \xi_3 (r) f.
    \end{cases}
    \]
    Choosing $q$ and $r$ as above, the system has no non-zero solutions and we have that $\kappa_{J_{q,r}} = -\infty$. \hfill $\blacksquare$
\end{example}

{\small
\printbibliography
}

\end{document}